\documentclass[a4paper,12pt,twoside]{article}

\usepackage{amssymb}
\usepackage[final]{lpretex}
\usepackage{title}
\usepackage{a4}
\usepackage{amscd}
\input xy
\xyoption{all}
\EndOfDump

\def\DHpartformat#1{(#1) }
\def\DHrefpart#1{(\DHRefpart{#1})}
\LBsetstyleof{partno}{arabic}

\DocVersion[Changed to LaTeX]{1}

\let\define\def

\def\A {{\mathbb A}}  \def\C {{\mathbb C}}
  \def\F {{\mathbb F}}

\def\GG {{\mathbb G}}   
  
  \def\P {{\mathbb P}} 
\def\Q {{\mathbb Q}}

\def\Z {{\mathbb Z}} 

\define \n {\mathbb N}
\define \z {\mathbb Z}
\define \q {\mathbb Q}
\define \PP {\mathbb P}

\def\sA {{\Cal A}}  
 \def\sE {{\Cal E}} \def\sF {{\Cal F}}
\def\sG {{\Cal G}} \def\sH {{\Cal H}} \def\sI {{\Cal I}}
\def\sJ {{\Cal J}}  \def\sL {{\Cal L}}
\def\sM {{\Cal M}} \def\sN {{\Cal N}} \def\sO {{\Cal O}}
 
 \def\sT {{\Cal T}} \def\sU {{\Cal U}}
  \def\sX {{\Cal X}}
\def\sY {{\Cal Y}}

\define \cN {\Cal N}
\define \cf {\Cal F}
\define \cg {\Cal G}
\define \cE {\Cal E}
\define \ce {\Cal E}
\define \cc {\Cal C}
\define \cV {\Cal V}
\define \cA {\Cal A}
\define \cK {\Cal K}
\define \cO {\Cal O}
\define \cF {\Cal F}
\define \cn {\Cal N}
\define \cI {\Cal I}
\define \sP {\Cal P}
\define \sW {\Cal W}
\def\tA {\widetilde{\Cal A}}

\def\a {\alpha} \def\b {\beta} \def\g {\gamma} \def \d {\delta} 
   
\def\s {\sigma}

\define \x {\xi}
\define \y {\eta}
\define \G {\Gamma}
\define \r {\rho}
\define \w {\omega}

\def\tX {\widetilde X}

\def \tV {\widetilde V}
\def \tI {\widetilde I}
\def \trho {\widetilde {\rho}}

\def \tp {\widetilde{\mathbb P}}
\define \tH {\widetilde H}
\define \tG {\widetilde{\Gamma}}
\define \tW {\widetilde W}
\define \tF {\widetilde F}
\define \tm {\widetilde m}
\define \St {\widetilde S}
\define \Xt {\widetilde X}
\define \tS {\widetilde S}
\define \tpsi {\widetilde \psi}
\define \tL {\widetilde L}
\define \tE {\widetilde E}
\define \tl {\widetilde l}
\define \tA {\widetilde A}
\define \tom {\widetilde\omega}
\define \tT {\widetilde T}
\define \tB {\widetilde B}
\define \tf {\widetilde f}
\define \tsA {\widetilde{\sA}}

\define \tsM {\widetilde{\sM}}
\define \tM {\widetilde M}
\define \tphi {\widetilde{\phi}}
\define \trho {\widetilde{\rho}}
\define \tR {\widetilde R}
\define \tp {\widetilde p}
\define \tq {\widetilde q}
\define \tc {\widetilde c}
\define \tsF {\widetilde {\sF}}
\define \tsN {\widetilde {\sN}}
\define \tsU {\widetilde {\sU}}
\define \th {\widetilde h}

\define\ti{\tilde\iota}

\def\pd {\partial}

\def \Dx1 {\frac{\pd}{{\pd} x_1}}
\def \Dy1 {\frac{\pd}{{\pd} y_1}}
\def \Dz1 {\frac{\pd}{{\pd} z_1}}
\def \Dx2 {\frac{\pd}{{\pd} x_2}}
\def \Dy2 {\frac{\pd}{{\pd} y_2}}
\def \Dz2 {\frac{\pd}{{\pd} z_2}}

\def\q {\quad} 

\def\mapdiagr#1{\Big\searrow\rlap{$\raise 5pt\vbox{{\hbox{$\mkern -15mu\scriptstyle#1$}}}$}}   

\def\mapdiagl#1{\llap{$\raise 5pt\vbox{{\hbox{$\scriptstyle#1\mkern
-15mu$}}}$}\Big\swarrow}              

\def\Mapdiagr#1{\nearrow\rlap{$\lower 5pt\vbox{{\hbox{$\mkern
-15mu\scriptstyle#1$}}}$}} 

\def\Mapdiagl#1{\llap{$\lower 5pt\vbox{{\hbox{$\scriptstyle#1\mkern
-15mu$}}}$}\searrow} 

\def\Mapswr#1{\swarrow\rlap{$\lower 5pt\vbox{{\hbox{$\mkern
-15mu\scriptstyle#1$}}}$}}              

\def\Mapnwl#1{\nwarrow\rlap{$\lower 5pt\vbox{{\hbox{$\mkern
-15mu\scriptstyle#1$}}}$}}

\def \inj {\hookrightarrow}

\define \Rhook {\hookrightarrow}

\def \half {\raise1pt\hbox{$\scriptstyle
        \frac{1}{2}\displaystyle$}}

\def \x{{\sl X}\llap{$\mkern -2mu {\scriptstyle -}$}}
\def \Proj {\operatorname{Proj}}
\def \Symm {\operatorname{Sym}}

\def \Pic {\operatorname{Pic}}
\def \Sing {\operatorname{Sing}}
\define \Kod {\operatorname{Kod}}
\define \dimension {\operatorname{dim}}
\define \codim {\operatorname{codim}}
\define \contr {\operatorname{contr}}
\define \rk {\operatorname{rank}}
\define \im {\operatorname{im}}
\define \Mor {\operatorname{Mor}}
\define \Cl {\operatorname{Cl}}
\define \Hilb {\operatorname{Hilb}}
\define \degree {\operatorname{deg}}
\define \mult {\operatorname{mult}}
\define \Aut {\operatorname{Aut}}
\define \NS {\operatorname{NS}}
\define \Gal {\operatorname{Gal}}
\define \ch {\operatorname{char}}
\define \Jac {\operatorname{Jac}}
\define \Km {\operatorname{Km}}
\define \Sec {\operatorname{Sec}}
\define \Stab {\operatorname{Stab}}
\define \Br {\operatorname{Br}}
\define \inv {\operatorname{inv}}
\define \tr {\operatorname{tr}}
\define \Frob {\operatorname{Frob}}
\define \Symn {\operatorname{Sym}^n}
\define \Ev {\sE^\vee}
\define \ordp {\operatorname{ord}_p}
\define \Supp {\operatorname{Supp}}
\define \Ann {\operatorname{Ann}}
\define \disc {\operatorname{disc}}
\define \Lie {\operatorname{Lie}}
\define \embdim {\operatorname{embdim}}

\def\id{\operatorname{id}}
\def\tj{\tilde{j}}
\def\tsJ{\widetilde{\sJ}}
\def\hV{\widehat{V}}
\def\bark{\overline{k}}

\def\barG{\overline{G}}
\def\barpi{\overline{\pi}}

\def\hod#1#2#3#4{\ensuremath{\if#30 H^{#2}({#1},{\cal O}_{#1}) \else 
 H^{#2}(#1,\Omega^{#3}\if\relax{#4}\relax_{#1}\else _{#1/#4}\fi)\fi}}

\begin{document}
\title
{Tropes, Torelli and theta characteristics}
\author{M. J. Fryers}
\address{School of Mathematics\\
University Walk\\
Bristol BS8 1TW\\
U.K.}
\email{michael@fryers.myzen.co.uk}
\author{N. I. Shepherd-Barron}
\address{Math. Dept.\\
King's College\\
Strand\\
London WC2R 2LS\\
U.K.}
\email{Nicholas.Shepherd-Barron@kcl.ac.uk}
\maketitle
\textup{2000} \textup{MSC}: 
14K25 (primary), 14H40 (secondary)
\begin{abstract}
This paper concerns the geometry associated to the $2$-torsion
subgroup $A[2]$ of a principally polarized abelian variety (ppav) $(A,\lambda)$.
The main results are these.

\begin{enumerate}
\item\label{(1)} (Theorem \ref{1.2}.)
If the characteristic is $2$ then
$(A,\lambda)$ is ordinary if and only if, for every symmetric theta divisor
$\Phi$ on $A$, there
is a $2$-torsion point that does not lie on $\Phi$.

\item\label{(2)} (Theorem \ref{serre} and Proposition \ref{recovery}.)
We give a geometrical description,
if either the characteristic of the ground field $k$ is zero or if
the genus of the curve is $3$,
of the quadratic twist observed by Serre \cite{LS}
that arises when $(A,\lambda)$ is geometrically isomorphic to a Jacobian.

\item\label{(3)} (Theorem \ref{main}.)
Suppose that $C$ is a non-hyperelliptic 
curve of genus $3$
in characteristic $\ne 2$ and that $Z$ is the set of its odd theta characteristics.
Then $Z$ is naturally embedded in $\P^6$ and the intersection
of quadrics though $Z$ is a normal del Pezzo surface $S$ of degree $2$.
The curve $C$ can then be recovered 
as the normalization of
the unique anti-sexcanonical curve 
on $S$ that is singular at each point of $Z$.
\end{enumerate}
\end{abstract}
\medskip
(\ref{(1)}) 
is due to Laszlo and Pauly \cite{LP}
when $(A,\lambda)$ is a Jacobian.

The geometrical description referred to in (\ref{(2)})
is in terms of a theorem of Welters \cite{W}
concerning the geometry of the linear system $\vert 2\Theta\vert$.
The hypothesis that the ground field be of characteristic
zero underlies much of the literature in this area
and we have not undertaken the task of trying to remove it.
In genus $3$ Beauville and Ritzenthaler \cite{BR} have
found a different geometrical description, under certain
additional hypotheses.

The index $r(S)$ of the del Pezzo surface in (\ref{(3)})
is $1$ if $\ch k\ne 3$
and $2$ if $\ch k=3$.
A del Pezzo surface of degree $2$ and index $1$
is anticanonically a double cover of $\P^2$
branched along a quartic; however, the quartic that appears here is not $C$
but rather the contravariant $K_1$ given in symbolical terms \cite{GY}
by $K_1=(abu)^4$ when the quartic $f$ defining $C$ is given
symbolically by $f=a_x^4=b_x^4$.

Lehavi \cite{L} has given another way of recovering $C$ from its bitangents;
our result appears to be complementary to his.

\begin{acknowledgments}
We are very grateful to Igor Dolgachev and Bjorn Poonen
for correspondence and discussions
on these matters.
\end{acknowledgments}

\begin{section}{Theta characteristics in characteristic two}
We start by recalling some basic results and notation.

Suppose first that $(Y,\Theta)$ is a principal symmetric abelian torsor (psat)
over a base $S$, with associated ppav
$(A,\lambda)$; recall (see \cite{SB}, for example)
that this means that $Y$ is projective over $S$,
that $\Theta$ is an effective ample Cartier divisor on $Y$, $Y$
is a torsor under $A:=\Aut^0_Y$,
that $\lambda:A\to\Pic^0_Y\cong\Pic^0_A$
defined by $\lambda(a)=t_a^*\Theta-\Theta$ is an isomorphism,
there is given an extension of the action
of $A$ on $Y$ to an action of the split extension
$A\rtimes[-1_A]$ and that
$\Theta$ is preserved by $\iota=[-1_A]$.
The standard example is $Y=\Pic^{g-1}_C$, $\Theta$
is the locus of effective classes in $Y$ and $A=\Pic^0_C$,
where $C$ is a curve of genus $g\ge 2$.

Of course, a point $y$ on $Y$ defines an isomorphism
$y':Y\to A$ given by subtraction. In particular, if $y\in T$,
then $y'$ takes $\Theta$ to a symmetric theta divisor
$\Phi_y$ on $A$, and all symmetric theta divisors on $A$
arise in this way.

Although, when $g\ge 2$, a psat is not naturally a ppav
and a ppav is not naturally a psat,
the natural
morphism from the stack $\sY_g$ of psat's $(Y,\Theta)$ to
the stack $\sA_g$
of ppav's $(A,\lambda)$ is an isomorphism. Then the
\emph{scheme of theta characteristics},
defined as the fixed locus $T:=Fix_Y=Fix_{[-1_A],Y}$
of $[-1_A]$ in $Y$,
is a torsor under the $2$-torsion $P:=A[2]$.
We also identify $T$ with the scheme that parametrizes
symmetric theta divisors $\Phi$ on $A$: given a point $t$ of $T$,
the corresponding divisor $\Phi_t$ is $\Phi_t=\Theta-t$.

Let $\sT\to\sY_g\cong\sA_g$
denote the universal scheme of theta characteristics.

The main result of this section is that, in characteristic $2$,
a ppav is ordinary if and only if on the corresponding
psat, there is a theta characteristic that does not lie
on the theta divisor.

The following result is merely a restatement
in a way that emphasizes psat's of
the discussion on pp. 132-135 of \cite{FC}.
Their stack $N_g$ of ppav's with a symmetric
theta divisor is isomorphic to the stack $\sT$
of psat's $Y$ with a point of $Fix_Y$.
(As usual, $e_n$ is the Weil pairing 
defined by $(Y,\Theta)$
on the $n$-torsion $A[n]$.)
\def\marf{\mathit{Arf}}
So fix a $g$-dimensional psat $(Y,\Theta)\to S$,
corresponding to the ppav $(A,\lambda)\to S$.
Put $P=A[2]$.  
Let $T'$ denote the scheme of morphisms
$t:P\to\mu_{2,S}$ such that
$$t(p+q)t(p)t(q)=e_2(p,q).$$
That is, the points of $T'$ are the $\mu_2$-valued quadratic
forms on $P$ whose polarization is $e_2$; \emph{cf.} [I], p. 214.
\begin{proposition}
\part $T'$ is naturally identified with $T$. 
\part There is a morphism $\marf:T\to\mu_2$
such that $\marf(t)\marf(t+p)=t(p).$
\part Over any geometric point $\s$ of $S$ of characteristic $\ne 2$,
$\marf(t)$ is the usual Arf invariant of $t$, regarded as a quadratic
form on the symplectic $\F_2$-vector space $P(\s)$.
\begin{proof}
We first construct the morphism $\marf$.

The fact that $\iota(\Theta)=\Theta$
means that we can identify $\sO_Y(-\Theta)$ with a
subsheaf, preserved by $\iota$, of the structure sheaf $\sO_Y$,
so with an $\langle\iota\rangle$-linearized invertible sheaf on $Y$.
So the sheaf $\sL=\sO_Y(\Theta)$ is $\langle\iota\rangle$-linearized.

Restrict to the subscheme $T$ of $Y$ and,
locally on $T$, pick a generator $s$ of $\sL\vert_T$.
Then $\iota^*(s)=u.s$ for some $u\in \sO_T^*$
and $s=\iota^*(u).\iota(s)=u^2.s$ since $\iota$ acts trivially on $T$.
So $u^2=1$.
\def \ts {\tilde{s}}
If $\ts$ is another local generator of $\sL\vert_T$, then
$\ts=w.s$ for some local section $w$ of $\sO_T^*$, and then
$$\iota^*\ts=\iota^*(w)\iota^*(s)=w.u.s=u.\ts$$
and therefore $u$ is defined as a morphism
$u:T\to\mu_{2,S}$; define $\marf=u$.
That is, for any local generator $s$ of $\sL\vert_T$
and for $t\in T$, 
$(\iota^*(s))(t)=\marf(t)s(t).$

Note that this construction commutes with base change.
 
(For any field $k$ of characteristic $\ne 2$
and any $k$-point $t$ of $T$, $\marf(t)=(-1)^{\mult_t(\Theta)}$.
This is proved in \cite{MEq}, bottom of p. $307$.)

Now define a morphism
$T\to\Mor(P,\mu_2)$ by $t(p)=\marf(t)\marf(t+p)$
for all points $t$ (not necessarily geometric) of $T$.
We need to show that $t(p)t(q)t(p+q)=e_2(p,q)$,
i.e., that $t$ is a $\mu_2$-valued quadratic 
form that returns the alternating form $e_2$.
To prove a formula such as this, nothing is lost by making a
faithfully flat base change $S'\to S$, so we can assume that
$T$ has an $S$-point $t_0$ and that therefore $Y=A$,
where $t_0$ is identified with $0_A$, and that
$\Theta$ is symmetric on $A$. This is the context of \cite{MEq}.

The $\langle\iota\rangle$-linearization of $\sL$ is 
an isomorphism $\phi:\sL\to\iota^*\sL$ covering $[-1_A]$.
Normalize $\phi$ by demanding that $\phi(0)=1$.
Then Mumford's morphism $e^\sL_*:P\to\mu_2$, defined by
$e^\sL_*(p)=\phi(p)$, is exactly
$e^\sL_*(p)=\marf(t_0).\marf(t_0+p)$.
Mumford proves, in all characteristics except $2$, that
$$e^\sL_*(p+q)=e^\sL_*(p)e^\sL_*(q)e_2(p,q).$$
So we can assume that $S$ has a closed point $\sigma$ of
characteristic $2$, and then that $(Y,\Theta)\to S$,
or $(A,\lambda)\to S$, is versal at $\sigma$ and that
$S$ is integral.
Then the formula to be proved holds over the generic point of $S$,
and so over all of $S$.

So the morphism $T\to\Mor(P,\mu_2)$ is a morphism
$T\to T'$. There is an action of $P$ on $T'$ given by
$(p(t')(q))=t'(q)e_2(p,q)$ for $t'\in T'$;
this makes $T'$ into a torsor
under $P$, and $T\to T'$ is then $P$-equivariant, so
an isomorphism.

The identification of $\marf$ with the usual Arf invariant
(in characteristic not $2$)
is proved, for Jacobians, in \cite{MTh}. From this, an argument involving
lifting to characteristic zero and the irreducibility of $\sA_g$
(only required here in characteristic zero) completes the proof.
\end{proof}
\end{proposition}
In all characteristics except $2$, it follows from
the properties of the Arf invariant, as is well known,
that there are just $2^{g-1}(2^g-1)$ points of $T$ with
odd multiplicity on $\Theta$ and $2^{g-1}(2^g+1)$ with
even multiplicity. That is,
provided that $\ch k\ne 2$,
$T=T^+\coprod T^-$ where $T^\pm=\marf^{-1}(\pm 1)$,
$T^+$ has order $2^{g-1}(2^g+1)$ and
$T^-$ has order $2^{g-1}(2^g-1)$.
 
The next result is due to Lazslo and Pauly \cite{LP}
in the case where $(A,\lambda)$ is a Jacobian.

\begin{theorem}\label{1.2} Assume that the base is
a field $k$ of characteristic $2$.

\part[i]\label{1.2(2)} There is at most one geometric 
point of $T$ that does not
lie on $\Theta$.

\part[ii] $(T\otimes\bark)_{red}$ is contained
in $\Theta\otimes\bark$ 
if and only if $A$ is not ordinary.

\part[iii] $A$ is ordinary if and only if, for every symmetric theta
divisor $\Phi_t$ on $A$, there is a geometric
point of $A[2]$ that does not lie on $\Phi_t$.
\begin{proof} 
We may assume that $k=\bark$.

Since $H^0(X,\sO(\Theta))$
is $1$-dimensional, any section $\theta$
satisfies
$\iota^*(\theta)=\pm\theta$. So in characteristic $2$,
$\iota^*(\theta)=\theta$. Recall that $A$ is ordinary
if and only if the identity connected component $P^0$ of $P$
has order at most $2^g$ (so exactly $2^g$).

Assume that $T_{red}$ does not lie in $\Theta$,
so that
there is a non-empty union $T^0=\sqcup T^0_i$
of connected components $T^0_i$ of $T$ that are disjoint from
$\Theta$. Then $P^0$ preserves each $T^0_i$
and $P^0$ has order $\ge 2^g$.
Also $\theta\vert_{T^0}$ generates $\sO(\Theta)\vert_{T^0}$,
and then, from the definition, 
$\marf =1$ on $T^0$. 

If $T^0$ is not connected, then there is a subgroup $P_1$ of $A[2]$
that contains $P^0$ as a subgroup of index $2$ and which preserves the union
$T^0_i\sqcup T^0_j$ of two distinct components $T^0_i$ and $T^0_j$.

Suppose that $R$ is a $k$-algebra,
that $t\in T^0(R)$ 
and that $p,q\in P^0(R)$. Then
$$e_2(p,q)=\marf(t)\marf(t+p)\marf(t+q)\marf(t+p+q);$$
each factor on the right equals $1$,
and so $P^0$ is totally isotropic.
Then $P^0$ has order $\le 2^g$
and $A$ is ordinary.

Moreover, if $T^0$ is not connected, then take $t\in (T^0_i\sqcup T^0_j)(R)$.
The same argument shows that $P_1$ is totally isotropic; this is impossible,
since its order is $2^{g+1}$.
So $T^0$ is connected, which proves \DHrefpart{i}.

Conversely, suppose that $A$ is ordinary.
It remains to prove that $T_{red}$ does not lie in $\Theta$.

Laszlo and Pauly \cite{LP} showed that on any ordinary
ppav $A$ over
an algebraically closed field of characteristic $p$, 
with symmetric theta divisor $\Phi=(\theta)_0$,
the translated powers $t_x^*\theta^p$ form a basis
of $H^0(A,\sO(p\Phi))$, where $x$ runs over the $k$-points
of $A[p]_{red}$. So if $\Phi$ contains $A[p]_{red}$, then
$A[p]_{red}$ is contained in the base locus of the linear
system $\vert p\Phi\vert$. But (Lefschetz) this base
locus is empty, and \DHrefpart{ii} is proved.

\DHrefpart{iii} follows at once from \DHrefpart{ii} 
and the correspondence between the symmetric theta divisors $\Phi_t$
on $A$ and the geometric points $t$ of $T$. 
\end{proof}
\end{theorem}
Let $\sA_g^{ord}$ denote the ordinary locus in
$\sA_g\otimes\F_2$ (regarded as the stack of psat's
in characteristic $2$) and let $\sT\to\sA_g^{ord}$
denote the universal scheme of theta characteristics.
(As mentioned above, this is isomorphic to the 
stack $N_g^{ord}$ of ordinary ppav's with
a symmetric theta divisor.)
Then $\sT=\sT^1\coprod\sT^2$,
where $\sT^1$ parametrizes theta characteristics
that do \emph{not} lie in the universal
$\Theta$ divisor and
$\sT^2$ parametrizes those that do.
So the restriction of $\marf$ to $\sT_1$
is identically $1$ and 
$\sT^1\to\sA_g^{ord}$ is a torsor
under the connected part of the
$2$-torsion subgroup of the universal abelian
variety.
On the other hand, the morphism
$\marf:\sT^2\to \mu_2$ is smooth
(\cite{FC}, p. 134).

The next lemma extends
Lemma 23 of [I], p. 214,
to characteristic $2$.
\begin{lemma}\label{igusa}
A decomposition $P=L\oplus M$
into Lagrangian subgroups 
determines a bilinear form $f:P\times P\to\mu_2$
such that $e_2(p,q)=f(p,q)f(q,p)$ and 
an even theta characteristic $\delta$
in $T(k)$ such that $\delta(p)=f(p,p)$.

For ordinary abelian varieties in
characteristic two, $\delta$
is the unique theta characteristic
that does not lie on $\Theta$.
\begin{proof}
We can identify $M=L^\vee$, the Cartier dual of $L$.
Define a bilinear morphism $f:P\times P\to\mu_2$
by $f((a,\a),(b,\b))=\b(a)$,
so that $e_2(p,q)=f(p,q)f(q,p)$,
and then define $\delta(p)=f(p,p)$.
So $\delta\in T(k)$.
Igusa's argument shows that
$arf(\delta)=1$ if $p\ne 2$,
while if $p=2$ then $arf(t)=1$
for any $k$-point $t$ of $T$.

If $A$ is ordinary in characteristic two,
then the local-{\'e}tale decomposition
of $A[2]$ is $A[2]=\mu_2^g\times(\Z/2)^g$,
and so determines $\delta$. The image of the
monodromy group on $A[2]$ is $GL_g(\Z/2)$,
so $\delta$ is the unique globally defined
theta characteristic.
\end{proof}
\end{lemma}
As showed to us by Bjorn Poonen,
the machinery of theta characteristics also illuminates
the finite Heisenberg subgroupschemes of a level $2$
theta group, as follows.

Recall that
if $(X,\Theta)$ is a psat
over some base $S$
and $(A,\lambda)$ is the corresponding ppav
then the \emph{theta group} $\sG_n=\sG_{(X,\Theta),n}$ 
of level $n$ attached to $(X,\Theta)$
is an $S$-groupscheme that is a central extension
$$1\to\GG_m\to\sG_n\stackrel{\pi}{\to} A[n]\to 0$$
from which the Weil pairing 
$e_n:A[n]\times A[n]\to\mu_n\inj\GG_m$ is then
constructed as the commutator pairing. The points of $\sG_n$
are pairs $(\phi,x)$ where $x\in A[n]$
and $\phi:t_x^*\sO_X(n\Theta)\to\sO_X(n\Theta)$ is an isomorphism
of line bundles. This description is given in \cite{MAV}
when $X$ is identified with $A$.

Part of the data of $(X,\Theta)$ is an involution
$\iota$ of $X$ that preserves $\Theta$
and is compatible with $[-1_A]$ acting on
$A=\Aut^0_X$.
This involution defines an involution
$\ti$ of $\sG_n$ given by
$\ti(\phi,x)=(\phi\circ\iota,-x)$.
The \emph{extended theta group}
$\sG_n^e$ is the split extension $\sG_n\rtimes\langle\ti\rangle$.

A \emph{finite Heisenberg group of level $n$ and type $(A[n],e_n)$} 
is a central extension
$$1\to\mu_n\to K\to A[n]\to 0$$
whose commutator pairing is $e_n$.

\begin{proposition} (Poonen)
The scheme $Heis_2=Heis_{(X,\Theta),2}$ that classifies 
finite Heisenberg subgroups
of $\sG_{(X,\Theta),2}$ that are of level $2$ and
type $(A[2],e_2)$
is a torsor over $S$ under $A[2]$.
It is isomorphic, as a torsor under $A[2]$, to the scheme
$T$ of theta characteristics.
\begin{proof}
Set $Z_2=\mu_2\subset\mu_4\subset\GG_m$
and $\sH_{4}=\{g\in\sG_2\vert g^{4}=1\}$.
It is easy to check that $\sH_4$ is a subgroup scheme
of $\sG_2$ and that $g^2\in\mu_2$ for all $g\in\sH_4$. 
Moreover, there is a central extension
$$1\to\mu_{4}\to\sH_{4}\to A[2]\to 0$$
and the finite Heisenberg subgroupschemes
of level $2$ and type $(A[2],e_2)$
of $\sG_2$ are exactly the subgroupschemes
$G$ of $\sH_{4}$ such that $G\cap\mu_{4}=Z_2$ and
$G\to A[2]$ is surjective. 
So $Heis_2$ is a locally closed subfunctor
of the Hilbert functor of $\sH_4$
and is therefore representable.
The extension
$$1\to\mu_{4}/Z_2\to\sH_{4}/Z_2\to A[2]\to 0 \ \ (*)$$
is an extension
of commutative groupschemes of type $(2,2,2,...)$
and, via replacing $G$ by
$\barG=G/Z_2$ and $\sH_{4}$ by $\overline{\sH_{4}}=\sH_{4}/Z_2$,
we see that $Heis_2$ is the scheme that parametrizes
the splittings of this last extension $(*)$
This exhibits the structure of $Heis_2$
as a pseudo-torsor under $\sH om(A[2],\mu_{4}/Z_2)$,
so, via 
the natural isomorphism
$\mu_{4}/Z_2\to\mu_2:s\mapsto s^2$
and the pairing
$e_2$, as a torsor 
under $A[2]$. 

To show that this pseudo-torsor is a torsor we can
assume that $S$ is the spectrum of an algebraically
closed field. We must find a splitting $\barG$
of the surjection $\overline{\sH_{4}}\to A[2]$.
This is equivalent to showing that the 
exact sequence
$$0\to A[2] \to\overline{\sH_{4}}^\vee \to \Z/2\Z\to0$$
that arises as
the Cartier dual
of the sequence $(*)$
is split. For this, just lift the element $1$ 
of $\Z/2\Z$ to any $S$-point of $\overline{\sH_{4}}^\vee$; 
that $S$-point is killed by $2$ since the whole group is.
The existence of this $S$-point shows that the pseudo-torsor is
a torsor.

We shall next
construct a morphism 
$\a:Heis_2\to T$.
\def\tx{\tilde{x}}
Take a subgroup $G$ of $\sH_{4}$
that is a point of $Heis_2$. Given $x\in A[2]$, 
choose $\tx\in G$ with $\pi(\tx)=x$; then $\tx^2\in\mu_2$. 
Note that if $\tx'\in G$ and $\pi(\tx')=x$,
then $\tx'=u\tx$ for some $u\in\mu_2$,
so that $\tx'^2=\tx^2$.
So there is
a morphism $t_G:A[2]\to\mu_2$ defined by
$t_G(x)=\tx^2$.

Verifying that $t_G(x+y)t_G(x)t_G(y)=e_2(x,y)$
for $x,y\in A[2]$ is immediate, and so there is a morphism
$\a:Heis_2\to T$ defined by $G\mapsto t_G$. We need to 
verify that
$\a$ is $A[2]$-equivariant.

Let $\s_G:G\to\barG$
and $\rho:\mu_{4}\to\mu_{4}/Z_2$ be the
quotients by $Z_2$,
and $\barpi_{\barG}:\barG\to A[2]$ the induced isomorphism.
Given a character $\chi:A[2]\to\mu_{4}/Z_2$,
put 
$$G_\chi=\{gs\vert g\in G,\ s\in\mu_{4},\ 
\rho(s)=\chi\circ\barpi_{\barG}\circ\s_G(g)\}.$$
Then $G_\chi\cap\mu_{4}=Z_2$, so that
$G_\chi$ is also a point of $Heis_2$.
Since $G_{\chi\psi}=(G_\chi)_\psi$ we
have an action of $A[2]$ on $Heis_2$.
This is the action referred to above.

Let $p\in A[2]$ and pick $g\in G$ such that
$\barpi_{\barG}\s_G(g)=p$.
Then
$t_G(p)=g^2$ and
$\barpi_{\barG_\chi}\s_{G_\chi}(gs)=p$
for some $s\in\mu_4$ with
$\rho(s)=\chi\barpi_{\barG}\s_G(g)$,
so that $t_{G_\chi}(p)=s^2g^2$.
Since the composite homomorphism
$\mu_4\stackrel{\rho}{\to}\mu_4/Z_2\stackrel{\cong}{\to}\mu_2$
equals the homomorphism $\mu_4\to\mu_2:s\mapsto s^2$,
it follows that $\chi\barpi_{\barG}\s_G(g)=s^2$.
So $\rho(s)=\chi(p)$ and $t_{G_\chi}=\chi t_G$,
as required.
\end{proof}
\end{proposition}

The morphism $\marf:Heis_2\to\mu_2$ induced by this isomorphism
$\a:Heis_2\to T$ is, over an algebraically closed
field of characteristic $\ne 2$, the morphism that distinguishes
between the two classes of extraspecial $2$-groups of order $2^{1+2g}$.
\end{section}
\begin{section}{Torelli's theorem and Serre's quadratic twist}
A crude version of the Torelli theorem for curves
states that, over an algebraically
closed field, a curve can be recovered from its 
Jacobian. A more precise version is given
by Oort and Steenbrink [OS].
We need some notation to state it.
They used the language of level structures
but we shall use that of stacks.

Recall that if $\sX$ is an algebraic stack with finite
stabilizers (for example, a Deligne--Mumford stack)
then $[\sX]$ denotes its coarse, or geometric, quotient.

On any pp abelian scheme $A\to S$, there is
an involution $[-1]$. So $\mathbb Z/2$ acts on
the moduli stack $\sA_g$. Stacks can be defined
as equivalence classes of groupoids; from this
point of view we 
define the quotient stack
$\tsA_g=\sA_g/(\Z/2)$ 
as follow. Recall that $\sA_g$ is the quotient
$X/R$ where $X$ is the disjoint union of
finitely many schemes of finite type over $\Sp\Z$,
$X$ is the base of an everywhere
versal family of pp abelian schemes
$(A,\lambda)\to X$, as in [FC],
and $R\to X\times_{\Sp\Z} X$ is the $Isom$
scheme $R=Isom_{X\times X}(pr_1^*(A,\lambda),pr_2^*(A,\lambda))$.
Then define $\tR=R/(-1)$ and
$\tsA_g=X/\tR$.
The quotient morphism
$\rho:\sA_g\to\tsA_g$ is a commutative 
gerbe banded by $\Z/2$,
so that, locally in the
{\'e}tale topology on $\tsA_g$, 
$\sA_g\cong\tsA_g\times B(\mathbb Z/2)$.

Equivalently, 
define the
prestack ${}_{pre}\tsA_g$ by $Ob({}_{pre}\tsA_g)=Ob(\sA_g)$
and $$Mor_{{}_{pre}\tsA_g}((A,\lambda),(B,\mu))=
Mor_{\sA_g}((A,\lambda),(B,\mu))/(-1).$$
Then $\tsA_g$ is the stack associated to
${}_{pre}\tsA_g$.

Assume that $g\ge 2$, and
consider the jacobian morphism $j_g:\sM_g\to\sA_g$.
This is given by sending a curve $C$ of genus $g$
to either the psat $(\Pic^{g-1}_C,\Theta)$
or the ppav $(\Pic^0_C,\lambda)$.

In a neighbourhood of a non-hyperelliptic curve,
$j_g$ is isomorphic to 
a morphism $X/G\to Y/(G\times (\mathbb Z/2)),$
of quotient stacks,
where $X$ and $Y$ are smooth, of dimensions $3g-3$
and $g(g+1)/2$ respectively,
while in  a neighbourhood of a hyperelliptic curve
$j_g$ is isomorphic to $X/H\to Y/H$ where
$\mathbb Z/2$ is a normal subgroup of $H$.
So there is a quotient stack $\pi:\sM_g\to\tsM_g$,
given locally by $X/G\stackrel{\id}{\to}X/G$
or $X/H\to [X/(\Z/2)]/(H/(\Z/2))$
such that $\tsM_g$ is normal, and is relatively normal over
$\Sp\mathbb Z[1/2]$, and $\pi$ is an isomorphism
on the non-hyperelliptic locus over $\Sp\mathbb Z$.
Moreover, there is a $2$-commutative diagram
$$
\xymatrix{
{\sM_g^{nh}}\ar@{^{(}->}[r]|{\circ}\ar[d]_{\cong}&{\sM_g}\ar[d]_{\pi}\ar[r]^{j_g}&
{\sA_g^{irred}}\ar[d]^{\rho}\ar@{^{(}->}[r]|{\circ}&{\sA_g}\ar[d]^{\rho}\\
{\tsM_g^{nh}}\ar@{^{(}->}[r]|{\circ}\ar[d]&{\tsM_g}\ar[r]^{\tj_g}\ar[d]&
{\tsA_g^{irred}}\ar@{^{(}->}[r]|{\circ}\ar[d]&{\tsA_g}\ar[d]\\
{M_g^{nh}}\ar@{^{(}->}[r]|{\circ}&{M_g}\ar[r]^{[j_g]}&
{A_g^{irred}}\ar@{^{(}->}[r]|{\circ}&{A_g}
}$$
where, as usual, $M_g=[\sM_g]$, $A_g=[\sA_g]$,
the superscripts $nh$ and $irred$
refer to non-hyperelliptic curves and
geometrically irreducible ppav's, respectively,
and 
$\xymatrix@1{\ar@{^{(}->}[r]|{\circ}&}$ 
denotes an open embedding.
(A psat $(Y,\Theta)\to S$ is irreducible if and only if
$\Theta$ is geometrically irreducible,
in the sense that every geometric fiber $\Theta_s$
is irreducible; this is equivalent
to the corresponding ppav $(A,\lambda)\to S$
being geometrically irreducible
as a ppav.)

Suppose that $\sX$ is a Deligne--Mumford stack over $\Sp\Z$, that
$\sY$ is a closed substack of $\sX$ and that $\sY$ is smooth
over $\Sp\Z$. Then we
say that $\sX$ has 
\emph{Veronese singularities}
along $\sY$ if locally in the {\'e}tale topology
there is 
\begin{enumerate}
\item a short exact sequence 
$$1\to \Z/2\to G\to H\to 1$$
of finite {\'e}tale groups;
\item
a smooth $\Z$-scheme $Y$;
\item
an equivariant action of $G$ on
$Y\times\A^N_\Z\to Y$ that preserves $Y\times\{0\}$
such that $\Z/2$ acts trivially on $Y\cong Y\times\{0\}$
and freely on $Y\times(\A^N_\Z-\{0\})$;
\item and an isomorphism
$\sX\to[(Y\times\A^N)/(\Z/2)]/H$
such that the composite morphism
$$\sY\to\sX\to [(Y\times\A^N)/(\Z/2)]/H\to Y/H$$
is an isomorphism.
\end{enumerate}

Let $\tsJ_g$ denote the image of $\tj_g$.
Then Oort and Steenbrink's version of the Torelli theorem [OS] is this.

\begin{theorem} (Torelli)
\part[i] The morphisms $j_g:\sM_g\to\sA_g^{irred}$ and 
$\tj_g:\sM_g\to\tsA_g^{irred}$ are finite and separate geometric points.
\part[ii] $\tj_g$ induces an isomorphism of automorphism group schemes.
\part[iii] $\tj_g$ is a closed embedding over $\Sp\Z[1/2]$.
\part[iv] $\tj_g$ is an embedding of the non-hyperelliptic locus
$\tsM_g^{nh}$.
\part[v] $j_g$ and $\tj_g$ induce closed embeddings
of $\sM_g^h$ and $\tsM_g^h$, respectively.
So $\pi$ induces an isomorphism of the
hyperelliptic loci.
\part[vi] $\tsM_g$ has Veronese singularities along the hyperelliptic locus.
\part[vii] At every closed point of $\tsJ_g^h$, the
Zariski tangent spaces of $\tsJ_g$ and $\tsA_g$
coincide.
\begin{proof}
Except for the statements concerning $\tsA_g$, this is 
nothing more than a translation of [OS].
The rest is then a simple observation.
\end{proof}
\end{theorem}

This has concrete corollaries, as follows.

\begin{corollary} Over $\C$, we can write
$\sA_g=\mathfrak H_g/Sp_{2g}(\Z)$
and $\tsA_g=\mathfrak H_g/PSp_{2g}(\Z)$,
where $\mathfrak H_g$ is the Siegel upper
half-space of degree $g$.
The locus of period matrices
in $\mathfrak H_g$ (that is, the inverse image
of $\tsJ_g$ in $\mathfrak H_g$) has Veronese
singularities along the hyperelliptic locus.
\noproof
\end{corollary}
\begin{corollary} (Serre) Suppose that $(A,\lambda)$
is a ppav over 
a field $k$ and that $K/k$
is a field extension for which there is a 
curve $C$ over $K$ such that
the Jacobian $JC$ is isomorphic to
$(A,\lambda)\otimes K$. Then there is a
curve $C_0$ over $k$ such that
$C$ is $K$-isomorphic to $C_0\otimes K$
and $(A,\lambda)$ is $k$-isomorphic to
the {\'e}tale quadratic twist of
$JC_0$ given by a unique quadratic character
$\epsilon$. If also $C$ is hyperelliptic,
then no twist is necessary.

Moreover, $C_0$ is unique up to a unique $k$-isomorphism.

Conversely, if $C$ is a non-hyperelliptic curve
and $(A,\lambda)$ is a non-trivial quadratic twist
of $JC$, then $(A,\lambda)$ is not a Jacobian.
\noproof
\end{corollary}
In fact this can be stated slightly more generally,
as follows.

\begin{corollary}\label{2.4}
Suppose that $S$ is a normal scheme and that
a family $(A,\lambda)\to S$
of irreducible ppav's over $S$ is given, 
defining $f:S\to\sA_g^{irred}$. 
Suppose also that there is a dense open subscheme
$i:S_0\inj S$
whose image in 
$A_g^{irred}$
lies in the image of $M_g^{nh}$.

Then there is a curve $C$ over $S$
such that $JC$ is $S$-isomorphic to an {\'e}tale
quadratic twist of $(A,\lambda)\to S$.
This twist is trivial along the inverse
image of the hyperelliptic locus in $S$.
\begin{proof}
We show first that there is a morphism $h:S\to\sM_g$
such that $\rho\circ f$ is isomorphic to $\rho\circ j_g\circ h$.

Since $M_g$ is the normalization of 
its image in $A_g^{irred}$, 
$S$ maps to $M_g$ in $A_g^{irred}$.
\begin{lemma}
$\tsM_g$ is identified with
the normalization of $(M_g\times_{A_g}\tsA_g^{irred})_{red}$.
\begin{proof}
Since $\tj_g$ is finite
the morphism
$\a:\tsM_g\to M_g\times_{A_g^{irred}}\tsA_g^{irred}=M_g\times_{A_g}\tsA_g$
is finite. 
Since $\a$ induces an isomorphism on geometric points and
$\tsM_g$ is normal
and $\tj_g$ induces an isomorphism of stabilizers,
the lemma follows.
\end{proof}
\end{lemma}
Since $S$ is normal $S\to M_g$ factors through $\tsM_g$,
say via $r:S\to\tsM_g$. 
Put 
$$\sY=\tsM_g\times_{\tsA_g}\sA_g=\tsM_g\times_{\tsA_g^{irred}}\sA_g^{irred}.$$
Then there is a $2$-commutative diagram
$$\xymatrix{
&S\ar[d]_{(r,f)}\ar@/^2pc/[dr]^f\ar@/^1pc/[dd]^>>>>>r&\\
{\sM_g}\ar@/_/[dr]_{\pi}\ar[r]_{(\pi,j_g)}\ar@/^3pc/[rr]^>>>>>>>>>>{j_g}
&{\sY}\ar[r]_{pr_2}\ar[d]_{pr_1}
&{\sA_g^{irred}}\ar[d]^{\rho}\\
&{\tsM_g}\ar[r]_{\tj_g}&{\tsA_g^{irred}.}
}$$
Now $S_0\inj S\stackrel{r}{\to}\tsM_g$ factors through
some morphism $h_0:S_0\to\sM_g^{nh},$
since $\sM_g^{nh}\to\tsM_g^{nh}$ is an isomorphism,
so there is a second $2$-commutative diagram
$$\xymatrix{
{S_0}\ar[dr]_{(h_0,i)}\ar@{^(->}[drr]|{\circ}^i\ar@/_2pc/[ddr]_{h_0}&&\\
&{\sM_g\times_{\sY}S}\ar[r]_<<<<{{\widetilde{pr}}_2}
\ar[d]^{{\widetilde{pr}}_1} &S\ar[d]^{(r,f)}\\
&{\sM_g}\ar[r]_{(\pi,j_g)}&{\sY.}
}$$
Since the composite
$j_g=pr_2\circ(\pi,j_g):\sM_g\to\sA_g^{irred}$
is finite, $(\pi,j_g)$ is also finite.
Therefore ${{\widetilde{pr}}_2}$ is finite,
and so, since $S$ is normal, 
$(h_0,i)$ extends to a morphism
$\d:S\to \sM_g\times_{\sY}S$.
That is, there is a morphism $h:S\to\sM_g$ that
extends $h_0$, and $\d=(h,1_S)$.
So the previous diagram can be extended to a third
$2$-commutative diagram
$$\xymatrix{
{S}\ar[dr]_{(h,1_S)}\ar[drr]^{1_S}\ar@/_2pc/[ddr]_{h}&&\\
&{\sM_g\times_{\sY}S}\ar[r]_<<<<{{\widetilde{pr}}_2}
\ar[d]^{{\widetilde{pr}}_1} &S\ar[d]^{(r,f)}\\
&{\sM_g}\ar[r]_{(\pi,j_g)}&{\sY.}
}$$
Comparison of the first and third diagrams shows that
\begin{eqnarray*}
\rho\circ f &\cong &\tj_g\circ pr_1\circ (r,f)
\cong\tj_g\circ pr_1\circ (r,f)\circ 1_S\\
{}&\cong &\tj_g\circ pr_1\circ (\pi,j_g)\circ h
\cong \rho\circ pr_2\circ (\pi,j_g)\circ h\cong\rho\circ j_g\circ h,
\end{eqnarray*}
as stated above.

Say $h(S)=(C\to S)$ and $(B,\mu)=JC$.
Put $I=Isom_S((A,\lambda),(B,\mu))$;
then $I\to S$ is finite and $(-1)$ acts freely on it.
An isomorphism $\rho\circ f\to\rho\circ j_g\circ h$
gives an $S$-point of $I/(-1)$.
Let $(X,\Theta)$ and $(Y,\Phi)$ be the psat's corresponding
to $((A,\lambda)$ and $(B,\mu)$, respectively;
then $I\cong Isom_S((X,\Theta),(Y,\Phi)$.
There are actions of $(-1)=\langle\iota\rangle$
on $(X,\Theta)$ and $(Y,\Phi)$ and $I$;
for $\g\in I$, $\iota(g)=[-1_B]\circ\g$.

Write $\tI=I/\langle\iota\rangle\to S$.
We know that $\tI$ has an $S$-point;
fix one such, say $S\to\tI$, and put $T=I\times_{\tI}S$.
This is a closed subscheme of $I$
and gives 
$$T\times_ST=\g_1\sqcup\g_2\inj I_T=Isom_T((X,\Theta)_T,(Y,\Phi)_T),$$
where $\g_{3-i}=[-1_B]\circ\g_i$.

Say $\Gal(T/S)=\langle\s\rangle$. Then $\s$ acts on
$(X,\Theta)\times_ST$ 
and $(Y,\Phi)\times_ST$ by $\s(x,t)=(x,\s(t)$,
$\s(y,t)=(y,\s(t))$.
So $\s\circ\g_i=\g_{3-i}\circ\s$.

Let $(Y',\Phi')$ denote the
quadratic twist of $(Y,\Phi)$ by $T\to S$.
Then $(Y',\Phi')\times_ST\cong (Y,\Phi)\times_ST$
and the action $\s*$ of $\s$ on 
$(Y,\Phi)\times_ST$ is given by
$\s*=[-1_B]\circ\s$. So
$$\s*\circ\g_i=[-1_B]\circ\s\circ\g_i=[-1_B]\circ\g_{3-i}\circ\s=\g_i\circ\s.$$
Therefore each $\g_i$ descends to an $S$-isomorphism
$\d_i:(X,\Theta)\to(Y',\Phi').$

This completes the proof of Corollary \ref{2.4}.
\end{proof}
\end{corollary}

A result of Welters \cite{W} gives a geometrical description
of Serre's quadratic character $\epsilon$, 
at least in characteristic zero or in genus three, as follows.
Welters' paper and many of those to which it refers
assume that the base field is $\C$; we have not 
checked whether this hypothesis is necessary.

Consider the subtraction map $s:C\times C\to\Jac^0_C:(P,Q)\mapsto [P-Q]$.
If $C$ is non-hyperelliptic, then $s$ is birational to its image $\Sigma$,
which has a unique singularity, the image of the diagonal $\Delta$. 

\begin{lemma} If $C$ is non-hyperelliptic, then $\Sigma$ is normal
and $C\times C$ is the blow-up of $\Sigma$ at the origin.

\begin{proof} Put $\frak m=\sI_{J,0}$. It is clear that
$\frak m.\sO_{C\times C}\subset \sI_\Delta$,
where $\sI_\Delta$ is the ideal sheaf of $\Delta$;
it is enough to show that this is an equality.

Since the natural map 
$$s^*\Omega_J^1\to\Omega^1_{C\times C}=pr_1^*\Omega^1_C\oplus
pr_2^*\Omega^1_C$$ 
induces $\omega\mapsto (\omega,-\omega)$ 
at the level of global sections, where $H^0(J,\Omega_J^1)$
is identified with $H^0(C,\Omega_C^1)$, it follows that
the homomorphism 
$\frak m/\frak m^2\to H^0(\Delta,\sI_\Delta/\sI_\Delta^2)=H^0(C,\Omega_C^1)$
is the identity. So
$\frak m.\sO_{C\times C}+\sI_\Delta^2=\sI_\Delta$,
and we are done, by Nakayama's lemma.
\end{proof}
\end{lemma}

In the natural 
$2\Theta$ linear system on $JC$,
let $\G_{00}$ be the linear subsystem consisting 
of those members that vanish to order
at least $4$ at the origin $0$. Consider
the intersection of the members of $\G_{00}$,
a subscheme of $JC$. Welters shows that,
up to embedded points and some $0$-dimensional material,
this subscheme is the surface $\Sigma$,
and $C$ is recovered as the exceptional divisor over $0$ in the minimal
desingularization $S\to\Sigma$. Moreover, the normality of $\Sigma$
means that $S$
is just the blow up of $\Sigma$ at the origin, so this
desingularization exists in families of non-hyperelliptic curves.

\begin{lemma} Suppose that $C$ is a curve of genus $g\ge 2$ and
that $S=C\times C$. Then the only morphisms
from $S$ to a curve whose generic fibre
has arithmetic genus $\g\le g$ are the two projections $pr_i:S\to C$.

\begin{proof} Suppose that 
the generic fibre $\phi$ of $q:S\to B$
has arithmetic genus $\g\le g$ and that 
$f_i$ is a fibre of $pr_i$ 
Then, by the adjunction formula,
$$2g-2=K_S.\phi+\phi^2 =(2g-2)(f_1+f_2).\phi,$$
so that, without loss of generality, $f_1.\phi=0$ and
$\phi$ is a fibre of $pr_1$, so that $q=pr_1$.
\end{proof}
\end{lemma}

\begin{theorem}\label{serre} In characteristic zero
Serre's quadratic character 
$\epsilon$ equals that given by the Galois action on the 
pair of projections $pr_i:S\to C$, $i=1,2$.
\begin{proof} It is clear that the involution $[-1_A]$ exchanges
the projections.
\end{proof}
\end{theorem}

We now verify this for curves of genus $3$ 
in all characteristics.

\begin{proposition} \label{recovery}
Suppose that $C$ is a non-hyperelliptic curve of
genus $3$ over a field
of any characteristic.
Then $\Sigma$ is the unique member
of $\vert 2\Theta\vert$ on $A=\Jac^0_C$ 
that vanishes to order at least $4$ at $0$. The curve $C$ is recovered
by resolving the singularity $(\Sigma,0)$
and $\epsilon$ has the same description as in Theorem \ref{serre}.

\begin{proof} First, it is well known that in characteristic zero,
$\Sigma\in\vert 2\Theta\vert$. Then the same result holds in characteristic
$p$, by specialization. Since $\mult_0\Sigma =4$, because the singularity
arises by contracting the plane quartic $\Delta$, $\Sigma$ lies in
$\Gamma_{00}$, the subspace of $2$nd order
theta functions vanishing to order
at least $4$ at $0$. It is enough to prove that
$\dim\Gamma_{00}=1$.

Suppose first that $p\ne 2$.
Let $\Gamma_0$ be the space of $2$nd order theta functions
vanishing at $0$; then $\dim\Gamma_0=7$ and every member of $\Gamma_0$
is even, so singular at $0$.
Since, in characteristic not $2$, 
the second order thetas provide an embedding
of the Kummer variety, the natural homomorphism
$\Gamma_0\to H^0(\P^2,\sO(2))$ is surjective. The kernel is
$\Gamma_{00}$, so that $\dim\Gamma_{00}=1$.

Now suppose that $p=2$. 
Since $\sO_A(2\Theta).\sO_{C\times C}=\sO_{C\times C}(pr_1^*K_C+pr_2^*K_C+2\Delta)$
and $\sI_{A,0}.\sO_{C\times C}=\sO_{C\times C}(-\Delta)$, it follows that
$$\sI_{A,0}^4\sO_A(2\Theta).\sO_{C\times C}=
\sO_{C\times C}(pr_1^*K_C+pr_2^*K_C-2\Delta),$$
so it's enough to show that 
$H^0(C\times C,\sO_{C\times C}(pr_1^*K_C+pr_2^*K_C-2\Delta))=0$.

There are exact sequences
$$0\to\sO(pr_1^*K_C+pr_2^*K_C-\Delta)\to\sO(pr_1^*K_C+pr_2^*K_C)
\to\sO_\Delta(2K_\Delta)\to 0,$$
$$0\to\sO(pr_1^*K_C+pr_2^*K_C-2\Delta)\to
\sO(pr_1^*K_C+pr_2^*K_C-\Delta)\to\sO_\Delta(3K_\Delta)\to 0.$$
Taking $H^0$ of the first sequence identifies
$H^0(C\times C,\sO(pr_1^*K_C+pr_2^*K_C-\Delta))$ with the kernel of
the natural multiplication
map $H^0(C,\sO(K_C))\otimes H^0(C,\sO(K_C))\to H^0(C,\sO(2K_C)),$
so with $\bigwedge{}^2H^0(C,\sO(K_C)).$ 
(Recall that, in any characteristic, the kernel of the natural projection
$\bigotimes^2V\to\Symm^2V$ is identified with
$\bigwedge{}^2V$ for any finite-dimensional vector space $V$.)

Taking $H^0$ of the second 
sequence then identifies the vector space
$$H^0(C\times C,\sO_{C\times C}(pr_1^*K_C+pr_2^*K_C-2\Delta))$$
with the kernel of the Wahl-Gauss map
$$\phi:\bigwedge{}^2H^0(C,\sO(K_C))\to H^0(C,\sO_C(3K_C)):\ 
s\wedge t\mapsto sdt -tds.$$ 

Note that, because
$h^0(C,\sO(K_C))=3,$ every element of
$\bigwedge{}^2H^0(C,\sO(K_C))$ is of the form $s\wedge t$.
Fix $0\ne\omega\in H^0(C,K_C)$, so that if $s,t\in H^0(C,K_C)$,
then $s=f\omega$ and $t=g\omega$ for rational functions
$f,g$ on $C$. Then $s\wedge t\in\ker\phi$ if and only if
$\frac{df}{f}=\frac{dg}{g}$.

Now $(f)=D-F$ and $(g)=E-F$, with $D,E,F\in\vert K_C\vert$
and $F=(\omega)_0$. Since $C$ is a double cover of its
Frobenius twist $C^{(1)}$, there are rational functions
$p,q,r,s,t$ on $C$ such that
$f=p^2+q^2t$ and $g=r^2+s^2t$; then
${df}/{f}={dg}/{g}$ if and only if
${r}/{s}={p}/{q}$. This implies that
${f}/{g}$ is a square; say ${f}/{g}=h^2$,
with $h\in k(C)$. Then there is an equality $2(h)=D-E$
of divisors on $C$; since $C$
is non-hyperelliptic, this leads at once to a contradiction.
\end{proof}
\end{proposition}

\begin{corollary} When $g=3$ the commutative diagram
$$\xymatrix{
{\sM_3^{nh}}\ar[r]\ar[d]_{\cong} & {\sA_3^{irred}}\ar[d]^{\rho}\\
{\tsJ_3^{nh}}\ar@{^{(}->}[r]|o & {\tsA_3^{irred}}
}$$
shows that the gerbe $\rho:\sA_3\to\tsA_3$ is neutral
over the open substack $\tsJ_3^{nh}$, so that
the stack $\sA_3^{smooth}$ of $3$-dimensional psat's
whose theta divisor is smooth is isomorphic to
$\sM_3^{nh}\times B(\Z/2)$. However, $\rho$
is not neutral over $\tsA_3^{irred}$
since $\sM_3\to\tsJ_3$ is not an isomorphism
in any neighbourhood of the hyperelliptic locus.
\noproof
\end{corollary}
\end{section}
\begin{section}{Tropes in genus $3$}
\emph{In this section the characteristic is not $2$.}

Suppose that $(Y,\Theta)$ 
is a psat of
dimension $g$ and that $(A,\lambda)$
is the corresponding ppav.

Assume also that $(Y,\Theta)$ is irreducible
as a psat (equivalently, that $(A,\lambda)$ is irreducible
as a ppav). Then
the complete linear system
$\vert 2\Theta\vert$ 
embeds the Kummer variety $\Km(Y):=[Y/[-1_A]]$
into $\P^{2^g-1}$.

There is a unique hyperplane $H$ in $\P^{2^g-1}$
such that $H.\Km(Y)=2X$, where $X=[\Theta/[-1_A]]$.
Both $H$ and $X$ are known as {\it{the trope}}
either of $Y$ or of $\Km(Y)$. Also, the 
$2$-torsion subgroupscheme $A[2]$ of $A$ acts
projectively on $\Km(Y)\subset \P^{2^g-1}$,
and every translate of $H$ or $X$ by a geometric point in $A[2]$
is also called a trope. However, we focus on psat's
rather than ppav's, and then there is a 
well defined choice of trope, as above.
The singular locus of $\Km(Y)$ is exactly the image
of $Fix_Y$ and the singular locus of $X$ is, if $g\ge 3$,
exactly the image of $\Sing\Theta\cup(\Theta\cap Fix_Y)$.

If $C$ is a curve of genus $g$ and $Y=\Jac^{g-1}_C$,
so that $A=\Jac^0_C$, then we write $\Km(C)$ instead of $\Km(A)$.

When $g=2$, that is, when $A$ is the Jacobian of a genus $2$ curve $C$, 
then the trope is a conic $\Gamma$ passing through $6$ of the $16$ nodes 
(ordinary double points, of type $A_1$) on
$\Km(A)$, and $C$ can be recovered, up to a quadratic twist, as the double cover
of $\Gamma$ ramified in the $6$ points \cite{Hu}.
Note that the ambiguity concerning $C$ in this quadratic twist is
the same ambiguity in recovering $(Y,\Theta)$ from the Kummer surface
and its trope. 
In fact, it is enough just to know the position of the $6$ nodes
in $\P^3$, since then $\Gamma$ is the unique conic through them;
it is the intersection of the quadrics in $\P^3$ through them.

We now show that in genus $3$ things are similar.

Recall that,
by definition, an index $1$ del Pezzo surface is a reduced and irreducible
Gorenstein surface
$S$ whose anti-canonical class $-K_S$ is ample. Its degree $\deg(S)$ is $K_S^2$.
Normal del Pezzo surfaces $S$
of index $1$ can be divided into two classes, as follows,
where $\tS\to S$ is the minimal resolution.

\noindent (1) Those with
only rational double points (RDPs). In this case 
either $\deg(S)=8$ and $S$ is a quadric or $1\le\deg(S)\le 9$ and 
$\tS$ is a blow-up of $\P^2$ in $9-\deg(S)$ points. 

\noindent (2) Those with a simply elliptic singularity. The degree may be 
any positive integer $d$
and $\tS$ is a $\P^1$-bundle over an elliptic curve $E$. 
The exceptional locus of $\tS\to S$ is a section $E_0$ of the bundle and $d=-E_0^2$;
$S$ can be regarded as the cone over an elliptic curve $E$ of degree $d$. 

If $S$ is a del Pezzo surface of index $1$ and degree $2$, then
$\vert-K_S\vert$ has no base points and defines
$S$ as a double cover of $\P^2$ branched in a quartic.
 

We shall also need to consider normal del Pezzo
surfaces of degree $2$ and of index $2$: 
for our purposes, these are defined as 
normal surfaces
$S$ that are quadric sections of the cone $\hV$ in $\P^6$
over a Veronese surface $V$ and that contain
the vertex $v$ of $\hV$.
For any such surface $-2K_S$ is a
very ample Cartier divisor, $K_S^2=2$,
and $S$ has a
single non-Gorenstein singularity, at $v$.

Recall the following result about tropes in genus $3$; see, for example,
Remark $6$ on p. 189 of \cite{DO}

\begin{proposition}\label{DO}
Suppose that $C$ is a non-hyperelliptic
curve of genus $3$ 
and that $X$ is a trope of $C$. 

\part[i] $X$ is a surface whose singular locus $Z$ consists of 
$28$ nodes. Its canonical class $K_X$ is ample,
$K_X^2=3$, $p_g(X)=3$, $q(X)=0$ and the complete linear system
$\vert K_X\vert$ has no base points. The canonical ring $R(X)$
is a hypersurface $k[x_1,x_2,x_3,y]/(f)$, 
where $\deg x_i=1$, $\deg y=2$ and $\deg f=6$.

\part[ii] The embedding $X\inj H\cong\P^6$ is defined by $\vert 2K_X\vert$.

\part[iii]\label{veronese} 
The intersection of the quadrics containing $X$ in $H$ is a 
copy $\hV$ of the cone over a Veronese surface $V$ in $\P^5$. 

\part[iv] The vertex $v$ of $\hV$ does not lie on $X$ 
and the morphism $X\to V\stackrel{\a}{\to} \P^2$ given by 
composing the projection from $v$ 
with an isomorphism $\a$
is defined by $\vert K_X\vert$.
\noproof
\end{proposition}
\begin{remark}
Horikawa showed \cite{Ho} that, in characteristic zero 
if $X$ is any canonical surface with $K_X^2=3$, $p_g(X)=3$, $H^1(X,\sO_X)=0$ 
and whose canonical system $\vert K_X\vert$ has no base points,
then $R(X)$ is as described. Ekedahl's 
results \cite{Ek} are enough to show that Horikawa's argument
applies to prove this in all characteristics.
\end{remark}

\begin{theorem}\label{main}
\part[i]\label{bitangent} 
The intersection in $H$ of the quadrics through $Z$ 
is a normal del Pezzo surface $S$ of degree 2, of index $1$ or $2$
and embedded in $H$ 
by $\vert -2K_S\vert$. 

\part[ii]
If $\ch k\ne3$, then $v\not\in S$ and $S$ is of index $1$.

\part[iii] If $\ch k=3$, then $v\in S$ and $S$ is of index $2$.

\part[iv] The curve $C$ is the normalization of 
the curve $R=X\cap S$. Moreover, $R$ is the unique curve
in $\vert-6K_S-2Z\vert$. That is, $R$ is the unique anti-sexcanonical curve lying
in $\vert-6K_S\vert$ and singular at $Z$.
\begin{proof}
Let $\pi:Y=\Pic^2_C\to\Km(C)$ be the quotient morphism
and consider the canonical morphism 
$\alpha=\phi_{\vert K_X\vert}:X\to \P^2$, which we identify
with projection from $v$ after fixing an isomorphism of the Veronese surface
$V$ with $\P^2$. According to Andreotti, $\alpha$ is separable
and $\deg\alpha=\frac{1}{2}{{2g-2}\choose{g-1}}=3$.

Suppose that $R\subset X$ is the ramification divisor and $B\subset\P^2$
the branch curve, the image of $R$, each given its reduced structure. 
Put $\tR=\pi^{-1}(R)$.
Then, by Andreotti's argument,
$\tR=R_1+R_2$ where
$R_1$ is the set of divisor classes $\sO(p+q)$ 
such that the line $\langle p,q\rangle$ is tangent to $C$
at a third point $r$, so that $K_C-p-q\sim 2r$,
and $R_2=\{\sO(2r)\vert r\in C\}$, the image of the diagonal.
The involution $\iota=[-1_A]$ exchanges $R_1$ and $R_2$,
so that each morphism $R_i\to R$ is 
an isomorphism outside the fixed locus
of $\iota$.

The doubling morphism $d:C\to R_2:r\mapsto \sO(2r)$ 
fits into a commutative diagram
$$\xymatrix{
{C}\ar[r]^d\ar[d]&{R_2}\ar[d]\\
{\Pic^1_C}\ar[r]^{[2]}&{\Pic^2_C}
}$$
where the vertical arrows are injective and $[2]$ is 
{\'e}tale,
and so separates tangent vectors. 
Since $d$ separates geometric
points ($C$ is not hyperelliptic)
$d:C\to R_2$ is an isomorphism.

Therefore the composite $\nu:C\to R_2\to R$ is 
an isomorphism \emph{except} that
$\nu(P)=\nu(Q)$ if $P+Q$ is a bitangent
and in this case $R$ has a node at $\nu(P)$,
and $R$ has a unibranched singularity at $\nu(P)$
if $P$ is a hyperflex of $C$.
So $R$ is smooth outside $Z$
and has points of multiplicity $2$ everywhere
on $Z$. Moreover,
$\nu$ is birational,
so that we shall be able to recover $C$ as the
normalization of $R$.

Note that, if $C$ is general,
then it has no hyperflexes, so that $R$ 
has only nodes and is smooth outside $Z$.
To see this, regard $C$ as the branch locus of a general
del Pezzo surface of degree $2$.

As Andreotti proved, it follows that
$B$ is the projective dual of the plane quartic $C$.

Recall that if $B$ is the dual curve of $C$
and $\tB$ is the normalization of $B$,
then $C\to\tB$ is some iterate of the Frobenius,
so that $p_g(B)=p_g(C)=3$ and
$$12=p^r\deg B$$
for some power $p^r$ of $\ch k$. 
So either $\deg B=12$, in which case
the morphisms $C\to R$ and $R\to B$ are both birational, 
or $\deg B=4$, $p=3$, $C\to R$ is birational and
$R\to B$ is birationally equivalent to the Frobenius,
so that $\deg(R\to B)=3$.
(If $p=2$ then $p^r\ge 2$, so that $\deg B=6$.)

\begin{lemma} If $C\to B$ is not birational,
then $p=3$ and $C$ is isomorphic to the Fermat
quartic.
\begin{proof}
Assume that $C\to B$ is not birational.
If $C$ is defined by $f=0$, then all the first partial
derivatives of $f$ are $p$th powers, so $p=3$
and the derivatives are cubes of linear forms.
Let $V_d$ denote the space of ternary $d$-ics,
regarded as a representation of $GL_3$.
Then $f$ lies in the
sub-representation $W=V_1\otimes V_1^{(1)}$
of $V_4$. Since $C$ is smooth,
the non-vanishing of the discriminant
and the finiteness of the automorphism
group show that $f$ is a semi-stable point
of $W$ whose stabilizer,
modulo the centre of $GL_3$, is finite.
So $f$ is a stable point. Since $\dim W=\dim GL_3$, there is
only one stable orbit.
\end{proof}
\end{lemma}

Moreover, the Pl{\"u}cker formulae give the following result.

\begin{lemma}
Suppose that $C$ has $\delta_1$ ordinary bitangents,
$\delta_2$ hyperflexes and $\kappa$ flexes.

\part[i] $\delta_1+\delta_2=28$.
\part[ii] The ordinary bitangents give nodes on $B$.
\part[iii] If $\ch k\ne 3$ then $2\delta_2+\kappa =24$.
The hyperflexes give unibranched singularities
whose value semigroup is generated by $\{3,4\}$
and the flexes give unibranched singularities
whose value semigroup is generated by $\{2,3\}$.
\part[iv] If $\ch k=3$ and $C$ is not the Fermat quartic then
$\delta_2+\kappa=8$.
The hyperflexes give unibranched singularities
whose value semigroup is generated by $\{3,5\}$
and the flexes give unibranched singularities
whose value semigroup is generated by $\{3,4\}$.
In particular, the non-nodal singularities of $B$
are of multiplicity $3$.
\begin{proof} \DHrefpart{i} is the basic
result about odd theta characteristics
and the rest is
based on the fact that, if $v(t)=(X:Y:Z)$
is a local parametrization of $C$, then
$v(t)\wedge v'(t)$ is a parametrization
of $B$.
\end{proof}
\end{lemma}

Let $\ell$ be the class of a line in $\P^2$. So
$$\alpha^*\ell\sim K_X\sim\alpha^*K_{\P^2}+e'R$$
for some $e'\ge 1$,
so that $4\alpha^*\ell\sim e'R$
and $e'(R.\alpha^*\ell)=12$.
Moreover, $R.\alpha^*\ell=\deg(R\to B).\deg B$,
so that $e'=1$ or $3$.

Suppose first that $e'=3$. Then $\deg B=4$ and $B$ is smooth,
and $R\to B$
is birational, so that $R$ is also smooth;
this is a contradiction.

Therefore $e'=1$.
Then $R$ is a quadric section of $X$ and $R.\alpha^*\ell =12$.
Since $R$ is singular at every point of $Z$ and $\deg Z=28$,
a manipulation of intersection numbers
on the minimal resolution $\tX$ of $X$ 
shows that $R$ is the unique
member of $\vert 4K_X-Z\vert$ (the point is that $(4K_X)^2< 28.2$). 
That is, $R$ is the unique quadric section
of $X$ passing through $Z$.

Regard $X$ as a Cartier divisor on the $3$-fold $\hV$ and let $\sO(1)$
denote the line bundle defining the embeddings of $X$ and $\hV$ into $H$.
So $\sO_X(1)=\a^*\sO_{\P^2}(2)$.
Then $H^0(\hV,\sO_{\hV}(2))\to H^0(X,\sO_X(2))$ is an isomorphism, so that
$R$ is cut out on $X$ by a unique member $S$ of $\vert\sO_{\hV}(2)\vert$.
Then $S$ is the unique member of $\vert\sO_{\hV}(2)\vert$ that
passes through $Z$, and so is the intersection in $H$ of the quadrics
through $Z$.

The reducedness and
irreducibility of $S$ follow from the fact that $R=S.X$
is a reduced, irreducible and ample Cartier divisor on $S$.

\begin{lemma} $S$ is smooth along $R$.
\begin{proof}
$R$ is smooth outside $Z$ and is Cartier on $S$,
so $S$ is smooth along $R-Z$.
At points of $Z$
the multiplicity of $R$ is $2$
and the multiplicity of $X$ is $2$. So $S$ has 
multiplicity $1$ along $Z$.
\end{proof}
\end{lemma}

\begin{corollary} $S$ is normal.
\begin{proof} $S$ is Cartier on $\hV$, so Cohen--Macaulay,
so satisfies Serre's condition $(S_2)$
and is smooth along an ample Cartier divisor,
so satisfies $(R_1)$.
\end{proof}
\end{corollary}

\begin{lemma} 
\part[i] $S$ contains $v$ if and only if $\ch k=3$.
\part[ii] $S$ is a degree $2$ del Pezzo
surface embedded by $\vert-2K_S\vert$.
Its index is $2$ if $\ch k=3$ and $1$ if $\ch k\ne 3$.
\begin{proof}
Write $\hV=\Proj k[x_1,y_1,z_1,w_2]$,
where the suffices indicate the degrees of the variables.
Then $X$ is defined in $\hV$
by the vanishing of a homogeneous sextic polynomial
$$F=w^3+w^2f_2(x,y,z)+wg_4(x,y,z)+h_6(x,y,z).$$
It follows that $S$ is defined in $\hV$
by the vanishing of the homogeneous quartic
$$\partial F/\partial w=3w^2+2wf+g,$$
so that $v\in S$ if and only if $\ch k=3$.
\end{proof}
\end{lemma}

Finally, consider $R$ as a curve on $S$. Since $R=S.X$ and is singular at $Z$, 
it is a member of
$\vert\sO_S(3) -2Z\vert$; i.e., it is a member of $\vert-6K_S-2Z\vert$.
Since $(6K_S)^2< 4.28$, $R$ is the only member of this linear system.
We have already remarked that $C$ is the normalization of $R$.
\end{proof}
This concludes the proof of Theorem \ref{main}.
\end{theorem}

\begin{proposition} Assume that $\ch k\ne 3$.

\part[i] 
$C_1$ is
the unique quartic in the dual
plane passing through the $24$ cusps
of the dual curve $B$ of $C$.
\part[ii] 
In symbolic terms \cite{GY},
if $C=a_x^4=b_x^4=c_x^4$,
then $C_1$
is given by $C_1=(abu)^4$.
\begin{proof}
We have already observed that $C_1$ is the unique quartic 
passing through the $24$ cusps of $B$.

Define contravariants $K_1$ and $K_2$ by
$$K_1=(abu)^4,\ K_2=(abu)^2(acu)^2(bcu)^2;$$
then it is well known that $B=K_1^3-6K_2^2$,
Here is a proof of this formula.

Consider a binary quartic $\a_y^4=\b_y^4=\g_y^4$.
Its ring of invariants is generated by $S$ and $T$,
where $S=(\a\b)^4$
and $T=(\a\b)^2(\a\g)^2(\b\g)^2$,
and its discriminant is $S^3-6T^2$.
Therefore
(this is the Clebsch transfer principle)
$K_1$ is the locus of lines $L$
such that $S(L\cap C)=0$, $K_2$
is the locus such that $T(L\cap C)=0$
and, by definition,
$B$ is the locus such that
$L\cap C$ is singular;
it follows that $B=K_1^3-6K_2^2$.

In particular, it is clear from this formula
that $B$ has cusps at every point of $K_1\cap K_2$
and that the intersection $K_1\cap K_2$ is transverse
(otherwise $B$ would have a singularity worse
than a cusp, which is not the case).
So $K_1\cap C_1$ contains at least $24$ points,
while $K_1.C_1=16$. So $K_1=C_1$, as required.
\end{proof}
\end{proposition}

\begin{remark}
For example, if $\ch k\ne 2,3$ and $C$ is the curve
$$X_1^3X_2+X_2^4+X_3^4=0,$$
then $C_1$ consists of $4$ lines through the point
$(0,1,0)$ in the dual plane, so that
the del Pezzo surface $S$ has a simply elliptic
singularity. This shows that the del Pezzo surface
$S$ is not enough to determine $C$, as follows.

Suppose
that $C$ is a contravariant of $C_1$. 
Every invariant of $C$ is then an invariant of $C_1$,
so vanishes, since $C_1$
is unstable. So $C$ is also unstable,
contradiction.
\end{remark}

\begin{remark} The expanded form of the symbolical expression for $C_1$
is as follows.
For non-negative integers $d,j_1,...,j_r$, write $j=(j_1,...,j_r)$ and let
${d\choose j}$ denote the multinomial coefficient
${d\choose j}=\frac{d!}{j_1!...j_r!}$ if $\sum j_i=d$ 
and ${d\choose j}=0$ otherwise.
If $C=\sum{4\choose i}A_ix^i$, where $i=(i_1,i_2,i_3)$
and $x=(x_1,x_2,x_3)$,
then 
$$C_1=\sum_\ell{4\choose\ell}(-1)^{\ell_2+\ell_4+\ell_6}
A_{\ell_1+\ell_2,\ell_3+\ell_4,\ell_5+\ell_6}
A_{\ell_4+\ell_5,\ell_1+\ell_6,\ell_2+\ell_3}u^{\ell_3+\ell_6,\ell_2+\ell_5,\ell_1+\ell_4}$$
where $\ell=(\ell_1,...,\ell_6)$ and $u=(u_1,u_2,u_3)$.
\end{remark}

\begin{remark} As mentioned above, Lehavi \cite{L} has already given an 
explicit way of recovering a smooth plane quartic $C$ from its bitangents. 
The comparison between his result and the one given above is this:
suppose given the locus $Z$ of $28$ points in the trope $H=\P^6$. 
Lehavi's result says that, if there is given also the
vertex $v$ of the Veronese cone $\hV$ that contains $X$, then 
there is an explicit procedure for recovering $C$
from the image of $Z$ in the Veronese surface $V$
that arises under projection from $v$.
On the other hand, our result does not demand the knowledge of $v$,
and it is not clear how to determine $v$ from knowledge solely 
of the $28$ points in $\P^6$.
\end{remark}

Now suppose that $C$ is hyperelliptic.
Most of the details 
are similar to those in the non-hyperelliptic case.
However, they are slightly more intricate because the surfaces $X$ and $S$ 
now pass through the vertex $v$ of the Veronese cone $\hV$ and have non-Gorenstein
singularities there. Recall that a singularity $X$ of type $[-4]$ is
a normal surface singularity the exceptional locus of whose minimal resolution
$\pi:\tX\to X$ consists of one smooth rational curve $E$
with $E^2=-4$. Such a singularity is a quotient $X=\A^2/\mu_4$, so rational, $2K_X$
is Cartier and $2K_{\tX}\sim \pi^*(2K_X)-E$.

\begin{theorem} 
Suppose that $C$ is a hyperelliptic curve of genus $3$ and $X$ its trope.

\part[i] $X$ is a surface whose singular locus is a set $Z$ of $28$ nodes
and one further point $v$ at which it has a 
singularity of type $[-4]$.

\part[ii] Let $\pi:\tX \to X$ be the minimal resolution of 
the singular point $v$.
Then $K_{\tX}$ is ample, $\tX$ has just $28$ nodes, $K_{\tX}^2=2$,
$p_g(\tX)=3$, $q(\tX)=0$ and $\vert K_{\tX}\vert$ has no base points.
The canonical morphism $\phi=\phi_{\vert K_{\tX}\vert}:\tX\to\P^2$ 
is of degree $2$
and the branch locus $B$ in $\P^2$ is the union of $8$ lines $B_i$ that are
tangent to a common conic $D$. The bicanonical map of $\tX$
is identified with the projection of $X$ from $v$ to $V$.

\part[iii] The composite morphism $\tX\to X\inj H\cong \P^6$ is defined by
the complete linear system $\vert 2K_{\tX}+E\vert$.

\part[iv] $X$ lies in a copy $\hV$ of a cone, with vertex $v$, 
over a Veronese surface $V$. 

\part[v] The intersection of the quadrics containing $X$ is $\hV$
and $X$ is a cubic section of $\hV$. In particular,
$X$ is a Cartier divisor on $\hV$.

\part[vi] The intersection of the quadrics through $Z$ is a normal surface $S$
that passes through $v$ and has a singularity of type $[-4]$ there.
It is a del Pezzo surface of index $2$ and degree $2$.
It is obtained 
by first taking the blow-up $\tS\to\P^2$ of $\P^2$ in $8$ distinct 
points $P_1,\ldots,P_8$
that lie on a conic ${\overline D}$
and then contracting the strict transform $D$ of ${\overline D}$.
The rational map $\P^2\to\P^6$ is defined by the linear system 
$\vert 4L-\sum P_i\vert$ of quartics through the $P_i$.

\part[vii] 
The embedding $S\inj H\cong\P^6$ is defined by $\vert-2K_S\vert$.

\part[viii] $X\cap S$ is the union of $8$ twisted cubics 
${\overline R}_i$, each passing through $v$.
As curves on $S$, each ${\overline R}_i$ is the strict transform of the 
tangent line to
${\overline D}$ at $P_i$ under the birational map referred to in \DHrefpart{vi}.

\part[ix] The divisor $\sum{\overline R}_i$ is Cartier on $S$ and is the unique member
of the linear system $\vert 3H -2Z\vert$.

\part[x] The curve $C$ is recovered by blowing up $S$ at $v$ to get 
an exceptional curve $E$
and then taking the double cover of $E$ ramified at the 
points where $E$ meets the strict
transforms of the curves ${\overline R}_i$.
\begin{proof} By the Riemann--Kempf singularity theorem, $\Theta$ has a node
$P$ at the point in $A[2]$ corresponding to the half--canonical $g^1_2$ on $C$.
There are three possible quotients of a node $(\Theta,P)$ by $[-1]$,
namely, a smooth point, 
an RDP of type $A_3$ and a $[-4]$ singularity, and they are distinguished by
whether $[-1]$ acts on the Zariski tangent space as 
a diagonal matrix $(-1,1,1)$, $(-1,-1,1)$ or $(-1,-1,-1)$.  
Since $[-1]$ acts on the tangent space $T_PA$ as $(-1)$, 
it follows that the quotient
$(X,v)=(\Theta,P)/[-1]$ is a $[-4]$ singularity. This is not Gorenstein,
but since $[-1]$ is an involution $2K_X$ is Cartier. The remaining singularities
of $X$ are a set $Z$ of $28$ nodes, as before.

Let $\pi:\Theta\to X$ be the quotient.
Since $[-1]$ acts freely in codimension one, $\pi^*(2K_X)=2K_\Theta$.
Moreover, exactly as in the non-hyperelliptic case, the morphism
$\Theta\to H=\P^6$ induced
by $\vert 2\Theta\vert$ factors through $X$ and $X\to H$ is 
defined by $\vert 2K_X\vert$ and is an embedding.
Also, the Gauss map on $\Theta$ is again, tautologically, defined
by $\vert K_{\Theta}\vert$ and factors through $X$. Since $X$ has
just rational singularities, it follows that $p_g(\tX)=3$.

From the description of the singularities on $X$, we get
$2K_{\tX}\sim\pi^*(2K_X)-E$. Suppose that $X^\dagger\to \tX$
is the minimal resolution and that
$F$ is a $(-1)$-curve in $X^\dagger$. Then $K_{X^\dagger}.F=-1$
and $\pi^*(2K_X).F>0$, so that $F.E\ge 3$. However, contracting
$F$ gives a curve, the image of $E$, that cannot live on a surface with $p_g>0$.
Hence $\tX$ is minimal, with $p_g =3$ and $K^2=2$. The classification
of such surfaces combined with Andreotti's proof of Torelli in the 
hyperelliptic genus $3$
case gives \DHrefpart{i} - \DHrefpart{iii}.

We next verify that $X$ is a cubic section of $\hV$; this will complete the proof of
\DHrefpart{iv} to \DHrefpart{v}.

Let $\sO(1)$ denote the hyperplane bundle on $\P^6$.
Then $\sO_{\tX}(3)=\sO_{\tX}(6K_{\tX}+3E)$, so that, by Riemann--Roch,
$\chi(\tX,\sO_{\tX}(1))=49$.
Considering exact sequences of the form
$$0\to\sO_{\tX}(6K_{\tX}+(n-1)E)\to\sO_{\tX}(6K_{\tX}+nE)\to\sO_E(6K_{\tX}+nE)\to 0$$
shows that $H^1(\tX,\sO(6K_{\tX}+3E)=0$, so that
$h^0(X,\sO_X(3))=49$. Since $\hV\cong\P(1,1,1,2)$,
it follows that $h^0(\hV,\sO(3))=50$. So $X$ lies in a non-trivial cubic section of
$\hV$. Since $\deg(X)=12=3\deg(\hV)$, $X$ is a cubic section of $\hV$. 
In particular, it is a Cartier divisor there.

Suppose that $R\subset\tX$ is the ramification curve of $\alpha:\tX\to\P^2$.
Since $2R=\alpha^*B$, it follows that $R=\sum R_i$ with $2R_i=\alpha^*B_i$. 
Since $B$ is tangent to $D$, we have $\alpha^{-1}(D)=E+E'$, where $E$ and $E'$
are exchanged by $[-1]$. So $E^2=E'^2=-4$ and $E.E'=8$. Also, if $L$
is a line in $\P^2$, then $K_{\tX}\sim\alpha^*L$ and
$$2=\alpha^*B_i.\alpha^*L=2R_i.K_{\tX},$$
so that $R_i.K_{\tX}=1$. 

Since  $2R_i.(E+E')=\alpha^*B_i.\alpha^*D=4$ and $R_i.E=R_i.E'$, by symmetry under
$[-1]$, we get $R_i.E=1$. Hence $R_i$ maps to a twisted cubic
in $\P^6$. Also, $2R\sim\alpha^*(8L)$, so that
$R+2E\sim 2\pi^*H\vert_{\P^6}.$ Define ${\overline R}$ to be the image of $R$ in $X$;
then ${\overline R}\sim 2H\vert_X$. That is, ${\overline R}$ is a quadric section,
passing through $Z$ and $v$.

Suppose that $T$ is another quadric section of $X$ through $Z$, not
necessarily through $v$.
Then in $\tX$ there are two members of $\vert 2K_{\tX}+E-Z\vert$,
namely $R+2E$ and $\pi^*T$. Since $(2K_{\tX}+E)^2=12$ and smooth curves
meeting transversely at a node have local intersection number $1/2$ there, we
have ${\overline R}=T$. So there is a unique quadric section $S$ of $\hV$
such that $X.S=R$ in $\hV$.

Hence $v\in S$ and $S$ is the intersection in $\P^6$ of the quadrics through $Z$.

We show next that $S$ is irreducible. For this, let $\gamma:\tV\to\hV$ be the
blow-up at $v$, with exceptional divisor $V_0\cong V$. Let $\tX$, $\tS$ be 
the strict transforms of $X$, $S$. Then $\gamma^*X\sim\tX+E$ and $\tX.V_0=E$.
The strict transform $\tR_i$ of ${\overline R}_i$ lies on $X$ and on $S$,
so that $\tR_i\subset \tX\cap\tS$ and $\tR_i$ meets $E$.
Hence $E$ and $S$ have at least $8$ points in common.

We have $\gamma^*S=\tS+aV_0$ for some $a\ge 1$. Now $E.V_0=-4$, 
since $E\inj V_0\cong\P^2$ realizes $E$ as a conic, so 
$0=\tS.E-4a$. 

Suppose that $E$ is \emph{not} contained in $\tS$;
then $E.\tS\ge 8$, so that $a\ge 2$.
We know, by looking at $X$, that projection $\beta$ from $v$ maps
$R_i$ to a line in $V\cong\P^2$, so $\beta(S)$ contains at
least $8$ distinct curves.
Let $\phi$ denote a fibre of the projection
$\tV\to V$. Then $\phi.V_0=1$ and $\phi$ maps to a line in $\P^6$,
so that $\phi.\gamma^*S=2$.
So, if $a\ge 2$, then $\phi.\tS=0$. That is, $\beta(S)\ne V$.
Then $S$ is of degree $8$ and has at least $8$ components.
However, $\hV$ contains no planes, and so $a=1$ and
$\phi.\tS=1$. Moreover, $E$ is contained in $\tS$.
 
So if $S$ is reducible, then $S=S_1+S_2$ with $S_1$ irreducible,
$v\in S_2$, $S_1$ is a linear section of $\hV$ and $S_2$ is a cone.
So $S_1\cong V$ and $S_2$ is a quartic cone. Since $E$ lies in $\tS$,
$S_2$ is irreducible. Then $R_1,\ldots,R_8$
are twisted cubics lying in $S_2$, which is impossible on the quartic cone $S_2$.

Having proved irreducibility, we move on to normality.
Since $S$ is an irreducible quadric section of $\hV$ through $v$,
the projection $\beta:S-\to V$ is birational. Thus $S$ can only fail to be normal
along generators of the cone $\hV$.

We have $\gamma^*S\sim \tS+V_0$, $\gamma^*X\sim \tX+V_0$, $\tX.V_0=E$ and 
$\tS\cap V_0\supset E$. So $\tS.V_0=bE+F$ for
some $b\ge 1$ and some curve $F\subset V_0$. So
$\tX.(bE+F)=\tX.\tS.V_0=\tS.E$
and
$$\tX.\tS.V_0= (\gamma^*X-V_0).(\gamma^*S-V_0).V_0=V_0^3=4.$$
Now $\tX.E=-V_0.E=4$, so that $b=1$ and $F=0$.
Hence $\tS$ is normal along $E$, and so $S$ is normal in a neighbourhood
of $v$, and so everywhere.

Now the identification $\hV=\P(1,1,1,2)$ and the adjunction formula
show that $-2K_S\sim H\vert_S$. Projection from $v$
completes the proof of \DHrefpart{vi} and \DHrefpart{vii}.

The only thing left is to show that ${\overline R}$ is the unique member
of $\vert 3H-2Z\vert$ on $S$. This is a straightforward calculation of 
intersection numbers
on $\tS$, as in the non-hyperelliptic case.
\end{proof}
\end{theorem}

\begin{remark} Suppose that the non-hyperelliptic curve $C_\eta$ specializes 
to a smooth hyperelliptic
curve $C_0$ and consider the corresponding specializations of 
tropes $X_\eta$ to $X_0$
and degree $2$ del Pezzo surfaces $S_\eta$ to $S_0$. Then $X_0$ and 
$S_0$ have singularities
of type $[-4]$, which are not Gorenstein, but from the viewpoint of the 
moduli of algebraic
surfaces of general type this is a well known picture: $X_0$ has 
semi-log canonical singularities and the
family $\{X_t\}$ is a relative canonical model, so that the surface $X_0$ 
appears on the boundary
of the separated and compact moduli space for stable surfaces with $K^2=3$ and 
$\chi(\sO)=4$.
\end{remark}
\begin{remark}
Tacit in this discussion so far has been the assumption that we began
with the Jacobian of some curve, and then sought to recover
the curve. However, there is another point of view.

Suppose that we
begin with a curve $\Gamma$ of genus $3$ over a field $k$, construct its
Jacobian $J\Gamma=A$ and then take the quotient $A/L$
of $A$ by
a Lagrangian subgroup $L$ of $A[n]$ for some $n\ge 2$. 
Given an equation for $\Gamma$ and sufficient knowledge
of torsion points on $A$, it is possible to write down equations
for $A$ such that the Kummer variety of the
principally polarized
quotient $A/L$ can also be written down
as the image of $A/L$ under its $2\Theta$
linear system. 

Then, by using the results of this paper,
we can construct a genus $3$ curve $C$
such that $A/L$ is a quadratic twist of $JC$
(provided that $A/L$ is geometrically irreducible as a ppav).
An analysis of the singular surface $\Sigma$ in the
linear system $\vert 2\Theta\vert$ will determine the
corresponding quadratic character. However, over a finite field $k$
this can often be done more easily by point-counting: the 
trace of Frobenius on $H^1(C\otimes\bark,\Q_\ell)$
is $\pm$ the trace of Frobenius on $H^1(\G\otimes\bark,\Q_\ell)$, so if this 
trace is non-zero the triviality or otherwise of the quadratic character is 
determined by a comparison of the two traces.
\end{remark}
\end{section}
\bibliography{alggeom,ekedahl}

\providecommand{\bysame}{\leavevmode\hbox to3em{\hrulefill}\thinspace}
\begin{thebibliography}{EGAIII:2}


\bibitem[BR]{BR}
A.~Beauville and C.~Ritzenthaler, \emph{Jacobians among abelian threefolds: a geometric approach},
Math. Ann. \textbf{350} (2011), 793--799.

\bibitem[DO]{DO}
I.~Dolgachev and D.~Ortland, \emph{Point sets in projective spaces and
theta functions}, Ast{\'e}risque \textbf{165} (1988).

\bibitem[Ek]{Ek}
T.~Ekedahl, \emph{Canonical models of surfaces of 
general type in positive characteristic}, Publ. Math. IHES \textbf{67} (1988), 97--144

\bibitem[FC]{FC}
G.~Faltings and C.-L.~Chai, \emph{Degeneration of abelian varieties},
Springer (1990).\textbf{67} (1988), 97--144.

\bibitem[GY]{GY}
G.H.~Grace and A.~Young, \emph{Algebra of invariants}, Cambridge, 1903.

\bibitem[Ho]{Ho}
E.~Horikawa, \emph{Algebraic surfaces of general type with small 
$c_{1}^{2}$, II}
Invent. Math. \textbf{37} (1976), 121--155.

\bibitem[Hu]{Hu} R.W.H.T.~Hudson, \emph{On Kummer's quartic surface},
Cambridge, 1905.

\bibitem[LP]{LP}
Y.~Laszlo and C.~Pauly, \emph{The Frobenius map, rank 2 vector bundles                                
and Kummer's quartic surface in characteristic 2 and 3}, Adv. Math.
\textbf{185} (2004), 246--269.

\bibitem[LS]{LS}
K.~ Lauter and J.P.~Serre,
\emph{Geometric methods for improving the upper bounds on the number of
rational points on algebraic curves over finite fields}. J. Algebraic Geom. 10 
(2001), no. 1, 19-36.

\bibitem[L]{L}
D.~Lehavi, \emph{Any smooth plane quartic can be 
reconstructed from its bitangents}, Isr. J. Math. \textbf{146} (2005),
371--379.

\bibitem[MEq]{MEq}
D.~Mumford, \emph{On the equations defining abelian varieties, I},
Invent. Math. \textbf{1} (1966), 287--354.

\bibitem[MAV]{MAV}
\bysame, \emph{Abelian varieties}, Tata Institute of Fundamental Research 
Studies in Mathematics, vol.~5, Oxford University Press., London, 1970.

\bibitem[MTh]{MTh}
\bysame, \emph{Theta characteristics of an algebraic curve},
Ann. Sci. ENS \textbf{4} (1971), 181--192.

\bibitem[SB]{SB}
N.~Shepherd-Barron, \emph{Thomae's formulae for non-hyperelliptic curves and 
spinorial square roots of theta-constants on the moduli space of curves},
arXiv:0802.3014

\bibitem[W]{W}
G.~Welters, \emph{The surface C-C on Jacobi varieties and 
2nd order theta functions}, Acta Math. \textbf{157} (1986), 1--22.

\end{thebibliography}
\bibliographystyle{pretex}
\end{document}